\documentclass[12pt]{amsart}
\usepackage{latexsym,amssymb,amsmath}
\usepackage[all,ps,cmtip]{xy} 
\newtheorem{theo}{Theorem}
\newtheorem{coro}{Corollary}
\newtheorem{prop}{Proposition}
\newtheorem{lemm}{Lemma}

\begin{document}

\def\ot{\otimes}
\def\we{\wedge}
\def\wec{\wedge\cdots\wedge}
\def\op{\oplus}
\def\ra{\rightarrow}
\def\lra{\longrightarrow}
\def\fso{\mathfrak so}
\def\cO{\mathcal{O}}
\def\cS{\mathcal{S}}
\def\fsl{\mathfrak sl}
\def\PP{\mathbb P}\def\PP{\mathbb P}\def\ZZ{\mathbb Z}\def\CC{\mathbb C}
\def\RR{\mathbb R}\def\HH{\mathbb H}\def\OO{\mathbb O}
\def\smc{\cdots }
\title{Chern classes of tensor products}
\author[L. Manivel]{Laurent Manivel}
\address{Institut Fourier, UMR 5582, 
Universit\'e Grenoble I and CNRS,
BP 74, 38402 Saint-Martin d'H\`eres, France}
\email{{\tt Laurent.Manivel@ujf-grenoble.fr}}

\begin{abstract} 
We prove explicit formulas for Chern classes of tensor products of vector bundles,
with coefficients given by certain universal polynomials in the ranks of the two bundles. 
\end{abstract}

\maketitle

\section{Introduction}

Chern classes are ubiquitous in algebraic topology, differential
geometry \cite{chern} or algebraic geometry  \cite{gr, fulton}.  
They have nice formal properties like the Whitney sum formula,
expressing the total Chern class of the direct sum of two 
complex vector 
bundles as the product of the total Chern classes of the two
bundles. The situation is much more complicated for the other
universal operation on vector bundles given by the tensor 
product: the Chern character is of course well behaved with respect
to products, but computing the Chern classes of the tensor product 
of two vector bundles is often a painful task. 

In this note we express the total Chern class of a tensor
product in terms of the Schur classes of the two bundles. 
Recall that the Schur classes are certain universal polynomials
in the Chern classes. They are indexed by partitions $\lambda
=(\lambda_1\ge\cdots\ge\lambda_r)$, and the Giambelli formula
expresses them as determinants in the usual Chern classes:
$$s_\lambda(E)=\det \big( c_{\lambda^*_i-i+j}(E) \big)_{1\le i,j\le s},$$
with the convention that $c_k(E)=0$ for $k<0$. Here $\lambda^*$
denotes the conjugate partition of $\lambda$, defined by 
$\lambda^*_i=\# \{k, \lambda_k\ge j\}$, and $s$ can be any 
integer greater or equal to $\lambda_1$. The Schur classes 
form an integral additive basis of the universal algebra generated 
by Chern classes, in particular there must be a universal formula
of type 
$$c(E\otimes F)=\sum_{\lambda,\mu}P_{\lambda,\mu}(e,f)s_\lambda(E)s_\mu(F)$$
for vector bundles $E,F$ of respective ranks $e,f$, the coefficients
$P_{\lambda,\mu}(e,f)$ being integers. In fact, the splitting 
principle allows to translate this identity into an identity of 
symmetric functions in two sets of variables, of size $e$ and $f$
respectively. 

An expression of this type has already been given by A. Lascoux 
in \cite{lascoux} (see also \cite[Ex.5 p.67]{mcd}), the coefficients $P_{\lambda,\mu}(e,f)$ being 
expressed as determinants of binomial coefficients. Explicitly:
$$P_{\lambda,\mu}(e,f)=
\det \Big( \binom{f-\mu^*_{e+1-i}+e-i}{\lambda_j+e-j}\Big)_{1\le i,j\le e}.$$
Unfortunately, these determinants seem quite 
difficult to evaluate in practice. Moreover, their dependence
in $e$ and $f$ appears quite unclear, while one can easily 
convince oneself that this dependence must be polynomial. 
Building on the work of Okounkov and Olchanski on shifted
Schur functions \cite{oo}, we obtain a  
polynomial formula for $P_{\lambda,\mu}(e,f)$. This
formula is very explicit, except maybe that it involves 
Littlewood-Richardson coefficients. Fortunately, our 
understanding of these fundamental coefficients has greatly 
improved in the recent years. In particular, many very
nice algorithms are known that allow to compute them quite
efficiently. 


\section{The main result} 

For any two partitions $\lambda$ and $\mu$, consider the polynomial
\begin{equation}\label{def}
P_{\lambda,\mu}(e,f)=\sum_\nu c_{\lambda^*,\mu}^{\nu^*}
(e|\nu-\lambda)(f|\nu^*-\mu)/h(\nu).
\end{equation}
The notation is the following. The partition $\lambda^*$ is the {\it conjugate
partition} of $\lambda$: when partitions are represented as  Young
diagrams, the lengths of the lines of $\lambda^*$ are the lengths
of the columns of $\lambda$. The coefficient $c_{\lambda^*,\mu}^{\nu^*}$
is a {\it Littlewood-Richardson coefficient} \cite{mcd}. It can be non-zero only 
when $\nu^*\supset\lambda^*$, or equivalently $\nu\supset\lambda$,
and $\nu^*\supset\mu$. The integer $h(\nu)$ is the product of the 
{\it hook-lengths} of the partition $\nu$, where the hook-length of a box 
$\alpha=(i,j)$ in $\nu$ is $h(\alpha)=\nu_i+\nu^*_j-i-j+1$. Finally, for a 
partition $\rho$, we let 
\begin{equation}\label{bingen}
(e|\rho)=\prod_{\alpha\in\rho}(e+c(\alpha)),
\end{equation}
where $c(\alpha)=j-i$ is the {\it content} of the box $\alpha=(i,j)$.
This is the  {\it content polynomial} of \cite[Ex.11 p.15]{mcd}. 
In particular $(e|k)=e(e-1)\cdots (e-k+1)$. This definition 
extends to skew-partitions: if $\rho\supset\sigma$, we 
simply let $(e|\rho-\sigma)=(e|\rho)/(e|\sigma)=\prod_{\alpha\in\rho/\sigma}(e+c(\alpha))$.

\medskip\noindent {\it Examples}. 
Suppose that $\lambda=(\ell)$ and $\mu=(m)$ have only one non-zero part. 
Then $\lambda^*=(1^\ell)$ has all its non-zero parts equal to one. The 
Littlewood-Richardson coefficient $c_{\lambda^*,\mu}^\nu$ is non-zero 
only if $\nu=(m,1^\ell)$ or  $\nu=(m+1,1^{\ell-1})$, and in both cases
it is equal to one. We thus get 
$$
\begin{array}{rcl}
P_{(\ell),(m)}(e,f)& =& \frac{(f-1)\cdots (f-\ell)(e+\ell)(e-1)\cdots (e-m+1)}{\ell !(m-1)!(\ell+m)}
 +\frac{(f+m)(f-1)\cdots (f-\ell+1)(e-1)\cdots (e-m)}{(\ell-1) !m!(\ell+m)},
\end{array}$$
$$P_{(\ell),(m)}(e,f)=\binom{e-1}{m-1}\binom{f-1}{\ell-1}\frac{ef-\ell m}{\ell m}.$$
Suppose now that $\lambda=(1^\ell)$ and $\mu=(1^m)$ have no part bigger than one.
By the previous computation and the symmetry properties stated in Proposition \ref{sym}, 
we get that 
$$P_{(1^\ell),(1^m)}(e,f)=\binom{e+m-1}{m-1}\binom{f+\ell-1}{\ell-1}\frac{ef-\ell m}{\ell m}.$$

The mixed case is more complicated. Suppose that $\lambda=(\ell)$ and $\mu=(1^m)$.
Using the symmetry properties of our polynomials we may suppose that $\ell\ge m$.  
Then the Littlewood-Richardson coefficient $c_{\lambda^*,\mu}^{\nu^*}$ is non-zero 
only if $\nu=(\ell+m-n,n)$ for some $n$ such that $0\le n\le m$, in which case
it is equal to one. We deduce the following formula:
$$P_{(\ell),(1^m)}(e,f)=\sum_{n=0}^m\frac{\binom{e+n-2}{n}\binom{e+\ell+m-n-1}{m-n}
\binom{f+1}{n}\binom{f-m}{\ell-n}}{\binom{\ell+m-n+1}{n}\binom{\ell+m-2n}{m-n}}.$$

For a last example, suppose that $\lambda=\mu=(2,1)$. Then $\nu$ is one of the 
partitions $(4,2), (4,1,1), (3,3), (3,2,1), (3,1,1,1), (2,2,2), (2,2,1,1)$. 
The corresponding Littlewood-Richardson coefficients are one, except for $\nu=(3,2,1)$, 
for which it is two. We get 
$$\begin{array}{rcl}
P_{(2,1),(2,1)}(e,f) & = & \frac{e(e-2)(e-3)f(f+2)(f+3)+e(e+2)(e+3)f(f-2)(f-3)}{80} \\
 & & +\frac{(e-3)(e-2)(e+2)(f-2)(f+2)(f+3)+(e-2)(e+2)(e+3)(f-3)(f-2)(f+2)}{72} \\
 & & +\frac{e(e-1)(e-2)f(f+1)(f+2)+e(e+1)(e+2)f(f-1)(f-2)}{144}
+2\frac{e(e-2)(e+2)f(f-2)(f+2)}{45}, \\
P_{(2,1),(2,1)}(e,f)  & =& \frac{e(e^2-1)f(f^2-1)}{9}-e^2f^2+e^2+2ef+f^2-4.
\end{array}$$

\begin{theo} 
Let $E,F$ be two vector bundles  of respective ranks $e,f$. The
total Chern class of their tensor product is 
$$c(E\otimes F)=\sum_{\lambda,\mu}P_{\lambda,\mu}(e,f)s_\lambda(E)s_\mu(F).$$
 \end{theo}

\proof By the splitting principle (see e.g. \cite[Remark 3.2.3]{fulton}), 
we are reduced to proving an identity between
symmetric functions in two sets of variables $x_1,\ldots ,x_e$ and $y_1,\ldots , y_f$. 
In the right hand side of the main formula, $s_\lambda(E)$ and $s_\mu(F)$ must 
then be replaced by the Schur functions $s_\lambda(x_1,\ldots ,x_e)$ and 
$s_\mu(y_1,\ldots , y_f)$ in these variables. 

In order to compute the right hand side, we start as in \cite{lascoux} 
with the Cauchy formula:
\begin{equation}
\prod_{\substack{1\le i\le e, \\ 1\le j \le f}}(1+x_iy_j)=
\sum_{\lambda\subset e\times f}s_{\lambda}(x_1,\ldots ,x_e)
s_{\lambda^*}(y_1,\ldots , y_f), 
\end{equation}
Replacing formally each $x_i$ by $x_i^{-1}$ 
and multiplying by $(x_1\ldots x_e)^f$ yields
\begin{equation}
\prod_{\substack{1\le i\le e, \\ 1\le j \le f}}(x_i+y_j)=
\sum_{\lambda\subset e\times f}s_{\lambda}(x_1,\ldots ,x_e)
s_{e\times f-\tilde{\lambda}}(y_1,\ldots , y_f).
\end{equation}
Here the notation is the following: the sum is over all partitions $\lambda$
whose Young diagram fits into the rectangle $e\times f$, which means that 
$\lambda_1\le f$ and $\lambda^*_1\le e$. Moreover we have denoted by 
$e\times f-\tilde{\lambda}$ the partition $(e-\lambda^*_f,\ldots , e-\lambda^*_1)$.
We deduce a first formula for the total Chern class of a tensor product:
$$\prod_{\substack{1\le i\le e, \\ 1\le j \le f}}(1+x_i+y_j)=
\sum_{\lambda\subset e\times f}s_{\lambda}(x_1,\ldots ,x_e)
s_{e\times f-\tilde{\lambda}}(1+y_1,\ldots , 1+y_f).$$
 
Now we use the {\it binomial theorem} \cite[Theorem 5.1]{oo} to obtain
$$s_{e\times f-\tilde{\lambda}}(1+y_1,\ldots , 1+y_f)=\dim_{GL(f)}(e\times f-\tilde{\lambda})
\sum_\mu\frac{s_\mu^*(e\times f-\tilde{\lambda})}{(f|\mu)}s_\mu(y_1,\ldots ,y_f).$$
Hence the following expression for the coefficient $P_{\lambda,\mu}(e,f)$ 
of $s_\lambda(E)s_\mu(F)$ in $c(E\otimes F)$:
$$P_{\lambda,\mu}(e,f)=\dim_{GL(f)}(e\times f-\tilde{\lambda})
\frac{s_\mu^*(e\times f-\tilde{\lambda})}{(f|\mu)}.$$
As in \cite{oo}, we have denoted by $\dim_{GL(f)}(e\times f-\tilde{\lambda})$ 
the dimension of the Schur module $S_{e\times f-\tilde{\lambda}}\CC^f$. It is 
given by the formula \cite[Ex.4 p.45]{mcd}:
\begin{equation}\label{dimschur}
\dim_{GL(f)}(e\times f-\tilde{\lambda})=\frac{(f|e\times f-\tilde{\lambda})}{h(e\times f-\tilde{\lambda})}.
\end{equation}
On the other hand $s_\mu^*$ denotes the {\it shifted Schur function} 
introduced in \cite{oo}. Its evaluation on a partition can be expressed in 
terms of representations of symmetric groups. Indeed, \cite[Theorem 8.1]{oo}
yields
$$s_\mu^*(e\times f-\tilde{\lambda}) = \frac{\dim [(e\times f-\tilde{\lambda})/\mu]}{\dim [e\times f-\tilde{\lambda}]}
(|e\times f-\tilde{\lambda}|\; |\; |\mu|).$$
Here $[\rho]$ denotes the irreducible representation of the symmetric group 
$\mathcal{S}_{|\rho|}$ associated to the partition $\rho$. Its dimension is
given by the celebrated hook-length formula \cite[Ex.2 p.74]{mcd}: letting $|\rho|=
\rho_1+\cdots +\rho_r$, 
\begin{equation}\label{hook}
\dim [\rho]= \frac{|\rho|!}{h(\rho)}.
\end{equation}
On the other hand $(e\times f-\tilde{\lambda})/\mu$ is not a partition but only a 
skew-partition, therefore the corresponding representation of the symmetric
group is not irreducible and there is no generalization of the hook-length 
formula that would give its dimension.
Nevertheless, its decomposition into irreducible representations is
known to be given by Littlewood-Richardson coefficients \cite[Ex.7 p.117]{mcd}:
\begin{equation}\label{lr}
[(e\times f-\tilde{\lambda})/\mu]=\bigoplus_\rho c_{\rho,\mu}^{e\times f-\tilde{\lambda}}[\rho].
\end{equation}
For this Littlewood-Richardson coefficient to be non-zero, we need that 
$\rho$ be contained in $e\times f-\tilde{\lambda}$. We can therefore write it as 
$\rho=e\times f-\tilde{\nu}$ for some partition $\nu$ containing $\lambda$. 
The coefficient $c_{\rho,\mu}^{e\times f-\tilde{\lambda}}$
is, by definition, equal to the multiplicity of the Schur module 
$S_{e\times f-\tilde{\lambda}}\CC^f$ inside the tensor product $S_\rho\CC^f\otimes S_\mu\CC^f$.
By \cite[Lemma 1]{man}, it is also equal to the multiplicity of $S_{f\times e}\CC^f=(\det\CC^f)^e$
inside the triple tensor product $S_\rho\CC^f\otimes S_\mu\CC^f\otimes S_{\lambda^*}\CC^f$.
But then for the same reason, it is also equal to the multiplicity of $S_{\nu^*}\CC^f$
inside $S_\mu\CC^f\otimes S_{\lambda^*}\CC^f$. In other words, 
we have proved that 
$$c_{\rho,\mu}^{e\times f-\tilde{\lambda}}=c_{\lambda^*,\mu}^{\nu^*}.$$
Therefore we get from (\ref{lr}) the identity
$$\frac{\dim [(e\times f-\tilde{\lambda})/\mu]}{\dim [e\times f-\tilde{\lambda}]}=
\sum_{\nu \subset e\times f} 
c_{\lambda^*,\mu}^{\nu^*}\frac{\dim [e\times f-\tilde{\nu}]}{\dim [e\times f-\tilde{\lambda}]}.$$
Using the hook-length formula (\ref{hook})
for $\dim [e\times f-\tilde{\nu}]$ and $\dim [e\times f-\tilde{\lambda}]$,
we deduce that 
\begin{equation}\label{quasi}
P_{\lambda,\mu}(e,f)=\frac{(f|e\times f-\tilde{\lambda})}{(f|\mu)}\sum_{\nu\subset e\times f}  
\frac{c_{\lambda^*,\mu}^{\nu^*}}{h(e\times f-\tilde{\nu})}.
\end{equation}

\begin{lemm}\label{prod}
$(f|e\times f-\tilde{\lambda})(e|\lambda)=(f|e\times f)=h(e\times f).$
\end{lemm}

\proof The quotient 
$(f|e\times f)/(f|e\times f-\tilde{\lambda})$ is the product of the $f+c(\alpha)$
for $\alpha$ a box in $e\times f$ but not in $e\times f-\tilde{\lambda}$. Such a box
has coordinates $\alpha=(f-j+1,e-i+1)$ with $1\le j\le \lambda_i$, and 
$f+c(\alpha)=f+(e-i+1)-(f-j+1)=e+j-i=e+c(\beta)$, where $\beta$ is a box
in $\lambda$.  Hence $(f|e\times f)/(f|e\times f-\tilde{\lambda})=(e|\lambda)$. 
The next identity is clear. \qed 

\medskip This leads for our coefficient  $P_{\lambda,\mu}(e,f)$
to the following expression:
\begin{equation}
P_{\lambda,\mu}(e,f)=\frac{1}{(e|\lambda)(f|\mu)}\sum_{\nu \subset e\times f} 
c_{\lambda^*,\mu}^{\nu^*}\frac{h(e\times f)}{h(e\times f-\tilde{\nu})}.
\end{equation}
Our next task will be to evaluate the quotient $h(e\times f)/h(e\times f-\tilde{\nu})$. 
In order to do this we will divide the rectangle $e\times f$ into four 
sub-rectangles NO, NE, SO, SE, in such a way that NO$\,\cup$ NE is the set 
of boxes $\alpha=(i,j)$ with $i\le f-\nu_1$, while NO$\,\cup$ SO is the set 
of boxes $\alpha=(i,j)$ with $j\le e-\nu^*_1$. We will denote by $h_{NO}(e\times f-\tilde{\nu})$, 
and so on, the product of the hook-lengths of the boxes of $e\times f-\tilde{\nu}$ belonging to
the rectangle NO. 

\begin{lemm}
The quotient $h(e\times f)/h(e\times f-\tilde{\nu})$ is the product of the 
following four partial quotients:
$$\begin{array}{rcl}
h_{NO}(e\times f)/h_{NO}(e\times f-\tilde{\nu}) & = & 1, \\
h_{NE}(e\times f)/h_{NE}(e\times f-\tilde{\nu}) & = & (e|\nu)/(\nu_1^*|\nu), \\
h_{SO}(e\times f)/h_{SO}(e\times f-\tilde{\nu}) & = & (f|\nu^*)/(\nu_1|\nu^*), \\
h_{NO}(e\times f)/h_{NO}(e\times f-\tilde{\nu}) & = & h(\nu_1\times\nu_1^*)/h(\bar{\nu}), 
\end{array}$$
where $\bar{\nu}$ denotes the partition $\nu_1\times\nu_1^*-\tilde{\nu}$.
\end{lemm}

\proof Straightforward. \qed

\medskip
We deduce a polynomial expression for our coefficient $P_{\lambda,\mu}(e,f)$:
\begin{equation}
P_{\lambda,\mu}(e,f)=\sum_{\nu\subset e\times f} 
c_{\lambda^*,\mu}^{\nu^*}(e|\nu-\lambda)(f|\nu^*-\mu)
\frac{h(\nu_1\times\nu_1^*)}{(\nu_1^*|\nu)(\nu_1|\nu^*)h(\bar{\nu})}.
\end{equation}
Indeed, this expression is really polynomial in $e$ and $f$ 
since we can omit the condition that $\nu$ be contained
inside the rectangle $e\times f$. If it is not, that is for example, 
if  $\nu_1^*$ is bigger
than $e$, then the box $\alpha=(e+1,1)$ is contained in $\nu$ and has content 
$c(\alpha)=-e$, which implies that $(e|\nu-\lambda)=0$. \medskip

In order to complete the proof of Theorem 1, there just remains to establish 
the following combinatorial lemma:

\begin{lemm}\label{id}
For any partition $\nu$, 
$$(\nu_1^*|\nu)(\nu_1|\nu^*)h(\bar{\nu})=h(\nu_1\times\nu_1^*)h(\nu).$$
\end{lemm}

\proof As $SL(\nu_1^*)$-modules, the Schur modules $S_\nu\CC^{\nu_1^*}$
and  $S_{\bar{\nu}}\CC^{\nu_1^*}$ are dual one to each other. In 
particular they have the same dimension, which means that 
$$\frac{(\nu_1^*|\nu)}{h(\nu)}=\frac{(\nu_1^*|\bar{\nu})}{h(\bar{\nu})}.$$
What  remains to notice is the identity $(\nu_1^*|\bar{\nu})=
h(\nu_1\times\nu_1^*)/(\nu_1|\nu^*)$, which is equivalent to
 Lemma \ref{prod}.\qed

\medskip\noindent {\it Remark}. Each term in  Lemma \ref{id}
is defined as 
a certain product of integers, and it seems that each integer $p$ 
appears the same number of times in the left and right hand sides of 
the identity. What is the combinatorial explanation? 

\medskip There is also a dual version of Theorem 1. Recall that total Segre class
of a vector bundle $E$ is defined as the formal inverse to the Segre class. More 
precisely, if we define the polynomial total Chern class of $E$ as 
$$c_t(E)=\sum_{k\ge 0}t^kc_k(E)=\prod_{i=1}^e(1+tx_i),$$
where $x_1,\ldots ,x_e$ are the formal Chern roots, then the polynomial total 
Segre class of $E$ is
$$h_t(E)=\sum_{k\ge 0}t^kh_k(E)=\prod_{i=1}^e(1-tx_i)^{-1}.$$
The total Segre class is $h(E)=h_1(E)$. 

\begin{theo} 
Let $E,F$ be two vector bundles  of respective ranks $e,f$. The
total Segre class of their tensor product is 
$$h(E\otimes F)=\sum_{\lambda,\mu}(-1)^{|\lambda|}P_{\lambda,\mu^*}(e,-f)s_\lambda(E)s_\mu(F).$$
 \end{theo}

The coefficient $Q_{\lambda,\mu}(e,f)=(-1)^{|\lambda|}P_{\lambda,\mu^*}(e,-f)$ 
of $s_\lambda(E)s_\mu(F)$ in this formula is
\begin{equation}
Q_{\lambda,\mu}(e,f)=\sum_\nu c_{\lambda,\mu}^{\nu}
(e|\nu-\lambda)(f|\nu-\mu)/h(\nu),
\end{equation}
and is clearly symmetric. 

\proof 
A completely formal argument shows that Theorem 1 is also valid for 
formal bundles. Indeed, first observe that the identity  $c(E\otimes (G\oplus H))=
c(E\otimes G)/c(E\otimes H)$ implies that the polynomials 
$P_{\lambda,\mu}(e,f)$ verify the relations
\begin{equation}\label{addition}
\sum_\mu c_{\varphi\psi}^\mu P_{\lambda,\mu}(e,g+h)=
\sum_{\alpha,\beta}c_{\alpha\beta}^\lambda P_{\alpha,\varphi}(e,g)P_{\beta,\psi}(e,h).
\end{equation}
This is a straightforward consequence of the fact that Littlewood-Richardson
coefficients also govern the decomposition of Schur classes of direct sums
\cite[I, (5.9)]{mcd}:
\begin{equation}\label{schuradd}
s_\mu(G\oplus H)=\sum_{\varphi,\psi}c_{\varphi\psi}^\mu s_\varphi(G)s_\psi(H).
\end{equation}
Now suppose that the formal bundle $F=G-H$, of rank $f=g-h$, is the formal difference of two vector
bundles $G$, $H$ of ranks $g$, $h$. Here $f=g-h$ can be negative. 
Then $E\otimes F=E\otimes G-E\otimes H$, hence $c(E\otimes F)=c(E\otimes G)/c(E\otimes H)$.
Theorem 1 for $F=G-H$ is thus equivalent to the identity
$$\begin{array}{rcl}
\sum P_{\lambda,\mu}(e,f)s_\lambda(E)s_\mu(F) &= &\sum
P_{\alpha,\beta}(e,f-g)s_\alpha(E)s_\beta(F-G)P_{\gamma,\delta}(e,g)s_\gamma(E)s_\delta(G) \\
 &=& \sum
P_{\alpha,\beta}(e,-h)P_{\gamma,\delta}(e,g)c_{\alpha,\gamma}^\theta s_\theta(E)s_\beta(F-G)s_\delta(G).
\end{array}$$
But (\ref{addition}) being a polynomial identity, remains valid if we replace 
$h$ by $-h$, and therefore the previous identity can be rewritten as
$$
\sum P_{\lambda,\mu}(e,f)s_\lambda(E)s_\mu(F)=
\sum
P_{\epsilon,\eta}(e,g-h)c_{\beta,\delta}^\eta s_\epsilon(E)s_\beta(F-G)s_\delta(G),$$
which clearly holds true since (\ref{schuradd}) is also valid for formal bundles, 
meaning that 
$$\sum_{\beta,\delta} c_{\beta,\delta}^\eta s_\beta(F-G)s_\delta(G)=s_\eta(F).$$

There just remains to apply 
Theorem 1, instead of $F$, to the formal bundle $-F$, of rank $-f$. We have
$c_t(-F)=h_{-t}(F)$, and more generally $s_\mu(-F)=(-1)^{|\mu|}s_{\mu^*}(F)$.
Therefore $h(E\otimes F)=c_{-1}(E\otimes (-F))$ is given by 
$$\begin{array}{rcl}
h(E\otimes F) & = & \sum_{\lambda,\mu}(-1)^{|\lambda|+|\mu|}P_{\lambda,\mu}(e,-f)s_\lambda(E)s_\mu(-F) \\
 & = & \sum_{\lambda,\mu}(-1)^{|\lambda|}P_{\lambda,\mu}(e,-f)s_\lambda(E)s_{\mu^*}(F).
\end{array}$$
This conclude the proof.\qed 




\section{Properties}

\subsection{Symmetries}

\begin{prop}\label{sym}
$P_{\lambda,\mu}(e,f)$ is an integer valued polynomial of degree $|\mu|$ in $e$ 
and degree $|\lambda|$ in $f$, with the following symmetries:
$$P_{\lambda,\mu}(e,f)=P_{\mu,\lambda}(f,e)=(-1)^{|\lambda|+|\mu|}P_{\lambda^*,\mu^*}(-e,-f).$$
\end{prop}

\proof The first assertion is obvious. To prove the first symmetry property 
we just need to notice that $c_{\lambda^*,\mu}^{\nu^*}=c_{\lambda,\mu^*}^{\nu}$
and $h(\nu)=h(\nu^*)$. To prove the second one we observe that if $\alpha$ is
a box of $\nu-\lambda$, then the corresponding box $\alpha^*$ 
in the conjugate skew-partition
$\nu^*-\lambda^*$ has opposite content. This implies that $(e|\nu-\lambda)=
(-1)^{|\nu|-|\lambda|}(-e|\nu^*-\lambda^*)$, and the conclusion easily follows.
\qed

\subsection{Vanishing}

\begin{prop}\label{zeroes}
One has $P_{\lambda,\mu}(e,f)=0$ whenever
$\lambda^*_1\le e<\mu_1$ or $\mu^*_1\le f<\lambda_1$.
\end{prop}

\proof 
If $c_{\lambda^*,\mu}^{\nu^*}\ne 0$, the Littlewood-Richardson rule 
implies that the first column of $\nu$ has length at least equal to 
$\mu_1$. If $\lambda^*_1<\mu_1$, this implies  that the intersection 
of $\nu-\lambda$ with the first column contains the boxes which belong 
to the lines numbered from $\lambda^*_1+1$ to $\mu_1$. These boxes 
have content $-\lambda^*_1,\ldots, -\mu_1+1$, hence $(e|\nu-\lambda)$
is divisible by $(e-\lambda^*_1)\cdots (e-\mu_1+1)$. Hence the first half
of the claim, the second one following by symmetry. \qed 

\subsection{Recursion}
Consider two complex vector bundles $E$, $F$ of respective rank $e$, $f$ and apply
Theorem 1 to $E'=E\oplus\cO$ and $F$, where $\cO$ denotes the trivial line 
bundle. Then $E'$ and $E$ have the same Chern and Schur classes. Since $E'\otimes F
=E\otimes F\oplus F$, the Whitney sum formula gives $c(E'\otimes F)=c(E\otimes F)c(F)$. 
Hence the relation
$$P_{\lambda,\mu}(e+1,f)= \sum_{\mu\rightarrow\theta}P_{\lambda,\theta}(e,f),$$
where $\mu\rightarrow\theta$ means that $\theta$ can be obtained from $\mu$ 
by suppressing some vertical strip. We can rewrite this as 

\begin{prop}\label{rec}
The polynomials $P_{\lambda,\mu}(e,f)$ obey the following recursion 
rule:
$$P_{\lambda,\mu}(e+1,f)-P_{\lambda,\mu}(e,f)= \sum_{\substack{\mu\rightarrow\theta, 
\\ \mu\ne\theta}}P_{\lambda,\theta}(e,f).$$
\end{prop}

We can use the same idea to obtain more recursion formulas. Indeed, 
suppose that $E=M\oplus \cO^{e-m}$
and $F=P\oplus \cO^{f-p}$ for some vector bundles $M,P$ or rank $m\le e$ and $p\le f$, 
respectively. Then $E$ and $M$ have the same Chern and Schur classes, as well as 
$F$ and $P$. The relation 
$$c(E\otimes F)=c(M\otimes P)c(M)^{f-p}c(P)^{e-m}$$
implies the following recursion formula, which is explicitly polynomial in $e$ and $f$:
$$P_{\lambda,\mu}(e,f) = \sum_{\alpha,\beta}P_{\alpha,\beta}(m,p)\sum_{\sigma,\tau}
d^{\lambda}_{\alpha,\sigma}d^{\mu}_{\beta,\tau}
\frac{(e-m|\tau_1+\cdots +\tau_p)}{\tau_1!\cdots \tau_p!}
\frac{(f-p|\sigma_1+\cdots +\sigma_m)}{\sigma_1!\cdots \sigma_m!}
.$$
Here we have denoted by $d^{\lambda}_{\alpha,\sigma}$ the generalized 
Kostka coefficient defined as 
the multiplicity of $s_\lambda$ inside 
 the product $s_\alpha e_1^{\sigma_1}\ldots e_m^{\sigma_m}$. 

\subsection{Leading term}

\begin{prop}\label{coef}
The leading term of $P_{\lambda,\mu}(e,f)$ is $e^{|\mu|}f^{|\lambda|}/h(\lambda)h(\mu)$.
\end{prop}

\proof Consider the previous formula for $P_{\lambda,\mu}(e,f)$. The term corresponding 
to the quadruple $\alpha,\beta,\sigma,\tau$ has degree $|\tau|=\tau_1+\cdots +\tau_p$ in $e$
and $|\sigma|=\sigma_1+\cdots +\sigma_m$ in $f$. But for $d^{\lambda}_{\alpha,\sigma}$
and $d^{\mu}_{\beta,\tau}$ to be non-zero we must have the relations $|\lambda|=|\alpha|
+\sigma_1+\cdots +m\sigma_m$ and $|\mu|=|\beta|
+\tau_1+\cdots +p\tau_p$. Hence  $|\tau|$ and $|\sigma|$ will be maximal 
when $\alpha,\beta$ are empty and $\tau_u,\sigma_v=0$  for $u,v>1$. But then 
the coefficient $d^{\lambda}_{\alpha,\sigma}$ is just the Kostka number $K_{\lambda}$,
the number of standard tableaux of shape $\lambda$. This is also the dimension of 
$[\lambda]$, and we can conclude the proof
by applying the hook-length formula (\ref{hook}) once again. \qed


\medskip
Comparing with the definition of $P_{\lambda,\mu}$ we deduce the following intriguing formula. 

\begin{coro}
For any three partitions $\lambda,\mu,\nu$, let 
$h^{\lambda,\mu}_\nu=h(\lambda)h(\mu)/h(\nu)$. Then
$$\sum_\nu h^{\lambda,\mu}_\nu c_{\lambda,\mu}^\nu=1.$$
\end{coro}

Is there any combinatorial interpretation ?

\end{document}